\documentclass[12pt,leqno]{article}
\usepackage{amsfonts}
\usepackage{amsmath, amsthm, amsfonts, amssymb, color}
\usepackage{mathrsfs}
\usepackage{color}
\theoremstyle{definition}

\newcommand{\scr}[1]{\mathscr #1}
\definecolor{wco}{rgb}{0.5,0.2,0.3}

\numberwithin{equation}{section} \theoremstyle{remark}

\newcommand{\ua}{\uparrow}

\title{{\bf Distribution Dependent SDEs  for   Navier-Stokes Type Equations }\footnote{Supported in
 part by  NNSFC (11831014, 11921001).} }
\author{{\bf  Feng-Yu Wang   }\\
\footnotesize{ Center for Applied Mathematics, Tianjin University, Tianjin 300072, China}\\
  \footnotesize{  Department of Mathematics,
Swansea University, Bay Campus, SA1 8EN, United Kingdom}\\
\footnotesize{    wangfy@tju.edu.cn}}
\begin{document}
\allowdisplaybreaks
\def\R{\mathbb R}  \def\ff{\frac} \def\ss{\sqrt} \def\B{\mathbf
B}
\def\N{\mathbb N} \def\kk{\kappa} \def\m{{\bf m}}
\def\ee{\varepsilon}\def\ddd{D^*}
\def\dd{\delta} \def\DD{\Delta} \def\vv{\varepsilon} \def\rr{\rho}
\def\<{\langle} \def\>{\rangle}
  \def\nn{\nabla} \def\pp{\partial} \def\E{\mathbb E}
\def\d{\text{\rm{d}}} \def\bb{\beta} \def\aa{\alpha} \def\D{\scr D}
  \def\si{\sigma} \def\ess{\text{\rm{ess}}}\def\s{{\bf s}}
\def\beg{\begin} \def\beq{\begin{equation}}  \def\F{\scr F}
\def\Ric{\mathcal Ric} \def\Hess{\text{\rm{Hess}}}
\def\e{\text{\rm{e}}} \def\ua{\underline a} \def\OO{\Omega}  \def\oo{\omega}
 \def\tt{\tilde}\def\[{\lfloor} \def\]{\rfloor}
\def\cut{\text{\rm{cut}}} \def\P{\mathbb P} \def\ifn{I_n(f^{\bigotimes n})}
\def\C{\scr C}      \def\aaa{\mathbf{r}}     \def\r{r}
\def\gap{\text{\rm{gap}}} \def\prr{\pi_{{\bf m},\varrho}}  \def\r{\mathbf r}
\def\Z{\mathbb Z} \def\vrr{\varrho} \def\ll{\lambda}
\def\L{\scr L}\def\Tt{\tt} \def\TT{\tt}\def\II{\mathbb I}
\def\i{{\rm in}}\def\Sect{{\rm Sect}}  \def\H{\mathbb H}
\def\M{\mathbb M}\def\Q{\mathbb Q} \def\texto{\text{o}} \def\LL{\Lambda}
\def\Rank{{\rm Rank}} \def\B{\scr B} \def\i{{\rm i}} \def\HR{\hat{\R}^d}
\def\to{\rightarrow}\def\l{\ell}\def\iint{\int}\def\gg{\gamma}
\def\EE{\scr E} \def\W{\mathbb W}
\def\A{\scr A} \def\Lip{{\rm Lip}}\def\S{\mathbb S}
\def\BB{\scr B}\def\Ent{{\rm Ent}} \def\i{{\rm i}}\def\itparallel{{\it\parallel}}
\def\g{{\mathbf g}}\def\Sect{{\mathcal Sec}}\def\T{\mathbb T}\def\BB{{\bf B}}
\def\f{\mathbf f} \def\g{\mathbf g}\def\BL{{\bf L}}  \def\BG{{\mathbb G}}
\def\Bd{{D^E}} \def\BdP{D^E_\phi} \def\Bdd{{\bf \dd}} \def\Bs{{\bf s}} \def\GA{\scr A}
\def\Bg{{\bf g}}  \def\Bdd{\psi_B} \def\supp{{\rm supp}}\def\div{{\rm div}}
\def\ddiv{{\rm div}}\def\osc{{\bf osc}}\def\1{{\bf 1}}\def\BD{\mathbb D}\def\GG{\Gamma}
\def\H{{\bf H}}
\maketitle

\begin{abstract}  To characterize Navier-Stokes type equations where the Laplacian is extended to a singular second order differential operator,  we propose a class of    SDEs depending on the distribution in future.    The well-posedness and  regularity estimates are derived for these SDEs.
  \end{abstract} \noindent

 \vskip 2cm

\section{Introduction }

 Let $d\in \mathbb N$. Consider the following  incompressible Navier-Stokes equation on $E:=\R^d$ or $\R^d/\Z^d$:
 \beq\label{E1} \pp_t u_t = \kk \DD u_t -(u_t\cdot\nn) u_t -\nn \wp_t,\ \ t\in [0,T] \end{equation}
 with  $\nn\cdot u_t:=\sum_{i=1}^d\pp_i u_t^i=0,$
 where  $T>0$ is a fixed time, $$u:=(u^1,\cdots, u^d): [0,T]\times E\to \R^d,\ \ \wp: [0,T]\times E\to \R,$$  and $u_t\cdot\nn:= \sum_{i=1}^d u_t^i \pp_i$.  This equation describes  viscous incompressible fluids,
  where $u$ is  the velocity field of a fluid flow, $\wp$ is the pressure, and $\kk>0$ is the viscosity constant.

Besides  existing probabilistic characterizations on Navier-Stokes equations,   see \cite{Qian} and references therein, in this paper  we propose a new type stochastic differential equation (SDE) depending on   distributions in the future, such that the solution  of \eqref{E1} is explicitly given by the initial datum $u_0$ and the pressure $\wp$. By proving the well-posedness of the SDE, we derive  the  well-posedness of \eqref{E1} in $\C_b^n (n\ge 2)$
with given pressure (which is however a part of solution in Navier-Stokes equations), see \cite{22n} for an analytic characterization  on the pressure to ensure $\nn\cdot u_t=0$.

Indeed, we will prove a more general result for    the following   Navier-Stokes type equation on $E:=\R^d$ or $E:=\R^d/\Z^d$:
 \beq\label{NS} \pp_t u_t =  L_t u_t -(u_t\cdot\nn) u_t +V_t,\ \ t\in [0,T],\end{equation}
 where
 $$L_t:= {\rm tr}\{a_t\nn^2\}+b_t\cdot\nn$$
 and
 $$V,\ b: [0,T]\times E\to\R^d,\ \ a: [0,T]\times E\to \R^{d\otimes d}$$ are measurable, and $a_t(x)$ is positive definite for  $(t,x)\in [0,T]\times E$.

 To characterize \eqref{NS}, we consider the following   SDE  on $\R^d$ where differentials are in $s\in [t,T]$:
 \beq\label{S1}\beg{split}&\d X_{t,s}^x=\ss{2 a_{T-s}}(X_{t,s}^x)\d W_s\\
  &+\bigg\{b_{T-s}(X_{t,s}^x)- \bigg[\E u_0(X_{s,T}^y)+\E\int_{s}^TV_{T-r}(X_{s,r}^y)\d r\bigg]_{y=X_{t,s}^x}\bigg\} \d s,  \\
 &\qquad  \ t\in [0,T], s\in [t,T], X_{t,t}^x=x\in\R^d,\end{split}\end{equation} where   $(W_s)_{s\in [0,T]}$ is a $d$-dimensional Brownian motion on a complete filtered  probability space $(\OO,\F,\{\F_s\}_{s\in [0,T]},\P)$.
 When $E=\T^d:=\R^d/\Z^d$, by extending a  function $f$ from domain $E$  to domain $\R^d$ as 
 \beq\label{PR} f(x+k)= f(x),\ \ x\in [0,1)^d, k\in \Z^d,\end{equation}
   we also have  the SDE \eqref{S1} for the case $E=\T^d$.

 Regarding $s$ as the present time, the SDE \eqref{S1}  depends on the distribution of $(X_{s,r})_{r\in [s,T]}$ coming from the future. So, this  is a future distribution dependent equation, but is essentially different from McKean-Vlasov SDEs 
 which depend on the distribution at present rather than future.  We will use $X:=(X_{t,s}^x)_{0\le t\le s\le T, x\in E}$  to formulate the solution to \eqref{NS}.

Let $D_T:= \{(t,s): 0\le t\le s\le T\}$. We define the solution $X$ of \eqref{S1} as follows.

\beg{defn}\label{D1} A family $X:=(X_{t,s}^x)_{(t,s,x)\in D_T\times \R^d}$  of random variables on $\R^d$ is called a solution of \eqref{S1}, if   $X_{t,s}^x$ is $\F_s$-measurable for all $x\in \R^d$ and $0\le t\le s\le T$,    $\P$-a.s.  continuous in $(t,s,x)$,
 $$\E\int_t^T  \bigg\{\big\|a_{T-s}(X_{t,s}^x)\|+ \bigg|b_{T-s}(X_{t,s}^x)- \bigg[\E u_0(X_{s,T}^y)+\E\int_{s}^TV_{T-r}(X_{s,r}^y)\d r\bigg]_{y=X_{t,s}^x}\bigg|\bigg\}\d s<\infty$$ for $ (t,x)\in [0,T]\times \R^d,$  and  $\P$-a.s.
\beg{align*}& X_{t,s}^x= x+   \int_t^s\ss{2 a_{T-r}}(X_{t,r}^x)\d W_r\\
&+ \int_t^s \bigg\{b_{T-r}(X_{t,r}^x)- \bigg[\E u_0(X_{r,T}^y)+\E\int_{r}^TV_{T-r}(X_{r,\theta}^y)\d \theta\bigg]_{y=X_{t,r}^x}\bigg\}\d r,\ \ (t,s,x)\in D_T\times\R^d.\end{align*}  \end{defn}

 We will allow the operator $L_t$ to be singular, where the drift  contains a locally integrable term introduced in \cite{XXZZ} for singular  SDEs.
 For any $p,q>1$ and $0\le t<s$, we write $f\in \tt L_q^p(t,s)$ if   $f=(f_r(x))_{(r,x)\in [t,s]\times\R^d}$ is a measurable function on $[t,s]\times\R^d$ such that
 $$\|f\|_{\tt L_q^p(t,s)}:=\sup_{z\in\R^d} \bigg(\int_t^s\|f_r1_{B(z,1)}\|_{L^p}^q\d r\bigg)^{\ff 1 q}<\infty,$$
 where $B(z,1)$ is the unit ball at $z$, and   $\|\cdot\|_{L^p} $ is the $L^p$-norm for the Lebesgue measure.
   We denote $f\in \tt H_{q}^{2,p}(t,s)$ if $|f|+|\nn f|+\|\nn^2 f\|\in \tt L_q^p(t,s).$ When $(t,s)=(0,T)$ we simply denote
   $$\tt L_q^p= \tt L_q^p(0,T),\ \ \tt H_q^{2,p}=\tt H_q^{2,p}(0,T).$$
 We will take $(p,q)$ from the following class:
 $$\scr K:=\Big\{(p,q): p,q>2, \ \ff {d}p +\ff 2 q<1\Big\}.$$
 We now make the following assumption on the operator $L_t$.

  \beg{enumerate}\item[$(H)$]  Let     $b_t=b^{(0)}_t+ b_t^{(1)}$, and when $E=\T^d$ we extend
   $a_t, b_t^{(0)}$ and $b_t^{(1)}$ to $\R^d$ as in $\eqref{PR}$.
  \item[$(1)$]  $a$ is positive definite  with
  $$\|a\|_\infty+\|a^{-1}\|_\infty:=\sup_{(t,x)\in [0,T]\times E}\|a_t(x)\|+\sup_{(t,x)\in [0,T]\times E}\|a_t(x)^{-1}\| <\infty,$$
$$ \lim_{\vv\to 0} \sup_{|x-y|\le \vv, t\in [0,T]} \|a_t(x)-a_t(y)\|=0.$$
 \item[$(2)$]   There exist $l\in \mathbb N$, $\{(p_i,q_i)\}_{0\le i\le l}\subset \scr K$ and $0\le f_i\in \tt L_{q_i}^{p_i}, 0\le i\le l$, such that
$$|b^{(0)}|\le f_0,\ \ \ \ \|\nn a\|\le \sum_{i=1}^l f_i. $$
 \item[$(3)$]          $\|b^{(1)}(0)\|_\infty:=\sup_{(t,x)\in [0,T]}|b^{(1)}(0)|<\infty $, and
 \beq\label{LP} \|\nn b^{(1)}\|_\infty:=\sup_{t\in [0,T]}\sup_{x\ne y} \ff{|b_t^{(1)}(x)- b_t^{(1)}(y)|}{|x-y|}  <\infty. \end{equation}
\end{enumerate}
Under this assumption, we will prove the well-posedness of \eqref{S1} and solve \eqref{NS} in the class
$$\scr U(p_0,q_0):=\Big\{u: [0,T]\times E\to\R^d;\
\,\|u\|_\infty+\|\nn u\|_\infty+\|\nn^2 u\|_{\tt L_{q_0}^{p_0}}<\infty\Big\}.$$
Recall that $W^{1,\infty}(E;\R^d)$ is the space of all weakly differentiable functions $f: E\to\R^d$ with $\|f\|_\infty+\|\nn f\|_\infty<\infty.$

\beg{thm}\label{T1} Assume $(H)$. Let $u_0\in W^{1,\infty}(E;\R^d)$ and $ \int_0^T\|V_t\|_\infty^2\d t<\infty.$  Then the following assertions hold.
\beg{enumerate} \item[$(1)$] The SDE \eqref{S1} has a unique solution $X:=(X_{t,s}^x)_{(t,s,x)\in D_T\times \R^d}.$
\item[$(2)$] If $u$ solves \eqref{NS} and $u\in \scr U(p_0,q_0)$,   then
 \beq\label{SL} u_t(x)= \E\bigg[u_0(X_{T-t,T}^x)+\int_{T-t}^T V_{T-s}(X_{T-t,s}^x)\d s\bigg],\ \ (t,x)\in [0,T]\times E.\end{equation} Moreover, there exists a constant $c>0$ such that for any $i\in\{1,2\}$ and $j,j'\in \{0,1\},$
\beq\label{EST} \|\nn^i u_t\|_\infty\le c t^{-\ff{i-j}2} \|\nn^j u_0\|_\infty + c\int_{T-t}^T (s+t-T)^{-\ff{i-j'}2}\|\nn^{j'} V_{T-s}\|_\infty \d s,\ \ t\in (0,T].\end{equation}
 \item[$(3)$] If $b^{(1)}=0$ and $ u_0, V_t\in \C_b^2$ with $\int_0^T \|V_t\|_{\C_b^2} \d t<\infty$,    then $u$ given by $\eqref{SL}$   solves $\eqref{NS}$, and   $u$ is in the class $\scr U(p_0,q_0)$.

\end{enumerate}
\end{thm}

 In the next two sections, we prove assertions (1) and (2)-(3) of Theorem \ref{T1} respectively, where in Section 2
 the well-posedness is proved for a more general equation than \eqref{S1}. Finally, in Section 4 we apply Theorem \ref{T1} to the equation \eqref{E1}.

\section{Proof of Theorem \ref{T1}(1)}

Let $\scr P$ be the set of all probability measures on $\R^d$ equipped with the weak topology, let $\L_\xi$ be the distribution of a random variable $\xi$ on $\R^d$.
Let
$$\GG:=C(D_T\times\R^d;\scr P)$$ be the space of continuous maps
 from $D_T\times\R^d$ to $\scr P$. For any $\ll>0$, $\GG$    is a complete     space under the metric
$$\rr_\ll(\gg^1,\gg^2):=\sup_{(t,s,x)\in D_T \times\R^d} \e^{-\ll(T-t)}\|\gg^1_{t,s,x}-\gg^2_{t,s,x}\|_{var},\ \ \gg^1,\gg^2\in \GG,$$ where   $\|\cdot\|_{var}$ is the total variation norm defined by
$$\|\mu-\nu\|_{var}:=\sup_{|f|\le 1}|\mu(f)-\nu(f)|,\ \ \mu,\nu\in \scr P $$
for $\mu(f):=\int_{\R^d} f\d\mu.$ Note that the convergence in $\|\cdot\|_{var}$ is stronger than the weak convergence.

We consider the following more general equation than \eqref{S1}:
   \beq\label{S2}\beg{split}&\d X_{t,s}^x= \Big\{b_{T-s}^{(1)}(X_{t,s}^x)+ Z_{s}(X_{t,s}^x,\L_{X}) \Big\}\d s + \ss{2 a_{T-s}}(X_{t,s}^x)\d W_s,\\
 &\qquad  \ t\in [0,T], s\in [t,T], X_{t,t}^x=x\in\R^d,\end{split}\end{equation}
where $\L_{X}\in \GG$ is defined by $\{\L_{X}\}_{t,s,x}  :=\L_{X_{t,s}^x},$ and
$$Z: [0,T]\times \R^d\times \GG\to \R^d$$ is measurable.

It is easy to see that \eqref{S2} covers \eqref{S1} for
\beq\label{ZZ} \beg{split}&Z_t(x, \gg):= b_{T-t}^{(0)}(x)- \int_{\R^d}u_0(y) \gg_{t,T,x}(\d y)-\int_t^T\d s\int_{\R^d} V_{T-s}(y)\gg_{t,s,x}(\d y),\\
&\qquad \ \ (t,x,\gg)\in [0,T]\times \R^d\times \GG.\end{split}\end{equation}
The solution of \eqref{S2} is defined as in Definition \ref{D1} using $b_{T-s}^{(1)}(X_{t,s}^x)+ Z_{s}(X_{t,s}^x, \L_{X})$ replacing
$$b_{T-s}(X_{t,s}^x)- \bigg[\E u_0(X_{s,T}^y)+\E\int_{s}^TV_{T-r}(X_{s,r}^y)\d r\bigg]_{y=X_{t,s}^x}.$$
We make the following assumption.

\beg{enumerate}\item[$(A)$]  $b^{(1)}$ and $a$ satisfy $(H)$, and there exists $(p_0,q_0)\in \scr K$ and $f_0\in \tt L_{q_0}^{p_0}$ such that
$$|Z_t(x,\gg)|\le f_0(t,x),\ \ (t,x,\gg)\in [0,T]\times \R^d\times \GG.$$
Moreover, there exists $0\le g\in L^2([0,T])$ such that
$$ \sup_{x\in\R^d} |Z_t(x,\gg^1)-Z_t(x,\gg^2)|\le g_t \sup_{(s,x)\in [t,T]\times \R^d} \|\gg^1_{t,s,x}-\gg_{t,s,x}^2\|_{var},\ \ t\in [0,T], \gg^1,\gg^2\in\GG.$$
\end{enumerate}

When $\|u_0\|_\infty+\int_0^T \|V_t\|_\infty^2\d t<\infty$,    $(H)$ implies $(A)$ for $Z$ given by \eqref{ZZ}. So,  Theorem \ref{T1}(1) follows from the following result, which also includes regularity estimates on the solution.

\beg{thm}\label{T2} Assume $(A)$. Then the following assertions hold.
\beg{enumerate}\item[$(1)$] $\eqref{S2}$ has a unique solution, and the solution has  the flow property \beq\label{FL} X_{t,r}^x= X_{s,r}^{X_{t,s}^x},\ \ 0\le t\le s\le r\le T, \ x\in\R^d.\end{equation}
\item[$(2)$] For any
$j\ge 1$, $$\nn_v X_{t,s}^{x}:=\lim_{\vv\downarrow 0} \ff{X_{t,s}^{x+\vv v}-X_{t,s}^{x}}{\vv},\ \ s\in [t,T]$$
exists in $L^j(\OO\to C([t,T]; \R^d), \P),$ and there exists a constant $c(j)>0$   such that
\beq\label{A*1} \sup_{(t,x)\in [0,T]\times \R^d} \E\bigg[\sup_{s\in [t,T]} |\nn_v X_{t,s}^{x}|^j\bigg]\le c(j)|v|^j,\ \ v\in \R^d.\end{equation}
\item[$(3)$] For any $0\le t<s\le T,$   $v\in\R^d$ and $f\in \B_b(\R^d)$,
\beq\label{BS} \nn_v \big\{\E f(X_{t,s}^\cdot)\big\}(x)= \ff 1 {s-t}\E\bigg[f(X_{t,s}^x)\int_t^s \Big\<\big(\ss{2 a_{T-r}}\big)^{-1}(X_{t,r}^x) \nn_v X_{t,r}^x,\ \d W_r\Big\>\bigg].\end{equation}\end{enumerate}
\end{thm}

\beg{proof} (a) We first explain the idea of proof using fixed point theorem on $\GG$.
For any $\gg\in \GG,$ we consider the following classical SDE
\beq\label{S3} \beg{split} &\d X_{t,s}^{\gg,x}= \Big\{b_{T-s}^{(1)}(X_{t,s}^{\gg,x})+ Z_{s}(X_{t,s}^{\gg,x},\gg)\Big\} \d s + \ss{2 a_{T-s}}(X_{t,s}^{\gg,x})\d W_s,\\
 &\qquad  \ t\in [0,T], s\in [t,T], X_{t,t}^{\gg,x}=x\in\R^d.\end{split}\end{equation}
By \cite[Theorem 2.1]{W21e} for $[t,T]$ replacing $[0,T]$, see also \cite{XXZZ} for $b^{(1)}=0$, this SDE is well-posed, such that for any $j\ge 1$ and $v\in\R^d$, the directional derivative
$$\nn_v X_{t,s}^{\gg,x}:=\lim_{\vv\downarrow 0} \ff{X_{t,s}^{\gg,x+\vv v}-X_{t,s}^{\gg,x}}{\vv},\ \ s\in [t,T]$$
exists in $L^j(\OO\to C([t,T]; \R^d), \P),$ and there exists a constant $c(j)>0$   such that
\beq\label{A*10} \sup_{(t,x)\in [0,T]\times \R^d} \E\bigg[\sup_{s\in[t,T]} |\nn_v X_{t,s}^{\gg,x}|^j\bigg]\le c(j)|v|^j,\ \ v\in \R^d,\end{equation} and  for any $f\in \B_b(\R^d)$,
\beq\label{BS'} \nn_v \big\{\E f(X_{t,s}^{\gg,\cdot})\big\}(x)= \ff 1 {s-t}\E\bigg[f(X_{t,s}^{\gg,x})\int_t^s \Big\<\big(\ss{2 a_{T-r}}\big)^{-1}(X_{t,r}^{\gg,x}) \nn_v X_{t,r}^{\gg,x},\ \d W_r\Big\>\bigg].\end{equation}
 By the pathwise uniqueness of \eqref{S3}, the solution satisfies the flow property
\beq\label{FL0} X_{t,r}^{\gg,x}= X_{s,r}^{\gg,X_{t,s}^{\gg,x}},\ \ 0\le t\le s\le r\le T, \ x\in\R^d.\end{equation}
Moreover,
$$\Phi(\gg)_{t,s,x}:=\L_{X_{t,s}^{\gg,x}},\ \ (t,s,x)\in D_T\times \R^d$$
defines a map $\Phi: \GG\to\GG$. If $\Phi$ has a unique fixed point $\bar\gg\in\GG$, then $\eqref{S3} $ with $\gg=\bar \gg$ reduces to \eqref{S2},   the well-posedness of \eqref{S3} implies that of \eqref{S2}, and  the unique solution is given by
$$X_{t,s}^x =X_{t,s}^{\bar\gg,x}.$$ Then   \eqref{FL}, \eqref{A*1} and \eqref{BS}  follow  from \eqref{FL0},
  \eqref{A*10}  and \eqref{BS'} for $\gg=\bar\gg$ respectively.   Therefore, it remains to prove that $\Phi$ has a unique fixed point.

(b) By the fixed point theorem, we only need to find constants $\ll>0$ and $\dd\in (0,1)$ such that
\beq\label{A*3} \rr_\ll(\Phi(\gg^1),\Phi(\gg^2))\le \dd \rr_\ll(\gg^1,\gg^2),\ \ \gg^1,\gg^2\in \GG.\end{equation}
Below, we prove this estimate using Girsanov's theorem.

For $i=1,2$, consider the SDE
 \beg{align*} &\d X_{t,s}^{i,x}= \Big\{b_{T-s}^{(1)}(X_{t,s}^{i,x})+ Z_{s}(X_{t,s}^{i,x},\gg^i)\Big\} \d s + \ss{2 a_{T-s}}(X_{t,s}^{i,x})\d W_s,\\
 &\qquad  \ t\in [0,T], s\in [t,T], X_{t,t}^{i,x}=x\in\R^d.\end{align*}
 By the definition of $\Phi$, we have
 \beq\label{A*4}\Phi(\gg^i)_{t,s,x}=\L_{X_{t,s}^{i,x}},\ \ i=1,2,\ (t,s,x)\in D_T\times \R^d.\end{equation}
Let
$$\xi_s:= \big(\ss{2a_{T-s}}(X_{t,s}^{1,x})\big)^{-1}\big\{ Z_{s}(X_{t,s}^{1,x},\gg^1)-  Z_{s}(X_{t,s}^{1,x},\gg^2)\big\},\ \ s\in [t,T].$$ By $(A)$, there exists a constant $K >0$ such that
\beq\label{A*5}|\xi_s|\le K g_s \sup_{(r,x)\in [s,T]\times \R^d}\|\gg^1_{s,r,x}-\gg^2_{s,r,x}\|_{var}.\end{equation}
By Girsanov theorem,
$$\tt W_s:= W_s -\int_t^s \xi_r\d r,\ \ s\in [t,T]$$ is a Brownian motion under the weighted probability
$\d\Q_t:= R_t\d\P$, where
$$R_t:= \e^{\int_t^T \<\xi_s,\d W_s\>-\ff 1 2 \int_t^T|\xi_s|^2\d s}.$$
With this new Brownian motion, the SDE for $X^1$ becomes
$$\d X_{t,s}^{1,x}=  \Big\{b_{T-s}^{(1)}(X_{t,s}^{1,x})+ Z_{s}(X_{t,s}^{1,x},\gg^2)\Big\} \d s + \ss{2 a_{T-s}}(X_{t,s}^{1,x})\d \tt W_s,\ \ s\in [t,T].$$
By the (weak) uniqueness for the SDE with $i=2$, we derive
$$\L_{X_{t,s}^{1,x}|\Q_t}=\L_{X_{t,s}^{2,x}}=\Phi(\gg^2)_{t,s,x},$$
where $\L_{X_{t,s}^{1,x}|\Q_t}$ is the distribution of $X_{t,s}^{1,x}$ under $\Q_t$. Combining this with
\eqref{A*4},   we get
\beq\label{P1}\|\Phi(\gg^1)_{t,s,x}-\Phi(\gg^2)_{t,s,x}\|_{var}=\sup_{|f|\le 1} \big|\E[f(X_{t,s}^{1,x})- f(X_{t,s}^{1,x})R_t]\big|
 \le \E|R_t-1|.\end{equation} By Pinsker's   inequality and the definition of $R_t$, we obtain
\beq\label{P2} (\E|R_t-1|)^2 \le 2\E [R_t\log R_t]= 2 \E_{\Q_t}[\log R_t]=2\E_{\Q_t}\int_t^T |\xi_s|^2\d s,\end{equation}
 where $\E_{\Q_t}$ is the expectation under the probability $\Q_t$.
 Combining \eqref{P1} and \eqref{P2} with \eqref{A*5}, and using the definition of $\rr_\ll$, we arrive at
 \beg{align*} & \|\Phi(\gg^1)_{t,s,x}-\Phi(\gg^2)_{t,s,x}\|_{var}\le   \bigg(2K^2  \int_t^T g_s^2\sup_{(r,y)\in [s,T]\times \R^d} \|\gg_{s,r,y}^1-\gg^2_{s,r,y}\|_{var}^2\d s\bigg)^{\ff 1 2}\\
 &\le \rr_\ll(\gg^1,\gg^2) \bigg(2K^2  \int_t^T g_s^2\e^{2\ll (T-s)}\d s\bigg)^{\ff 1 2}, \ \ (t,x)\in [0,T]\times \R^d.\end{align*}
 Therefore
 $$\rr_{\ll}(\Phi(\gg^1),\Phi(\gg^2))\le \vv_\ll \rr_\ll(\gg^1,\gg^2),$$
 where
 $$\vv_\ll:= \sup_{t\in [0,T]} \bigg(2K^2  \int_t^T g_s^2\e^{-2\ll (s-t)}\d s\bigg)^{\ff 1 2}\downarrow 0\ \text{as}\ \ll\uparrow\infty.$$
 By taking large enough $\ll>0$,   we prove \eqref{A*3} for some $\dd<1$.

\end{proof}

For later use we present the following  consequence of Theorem \ref{T2}.

\beg{cor}\label{CN} Assume $(A)$ and let
$$P_{t,s}f(x):= \E[f(X_{t,s}^x)],\ \ (t,s,x)\in D_T\times \R^d, f\in \B_b(\R^d).$$
Then  there exists a constant $c>0$ such that for any function $f$,
\beg{align*} &\|\nn P_{t,s}f\|_{\infty}\le c \min\Big\{(s-t)^{-\ff 1 2 } \|f\|_{\infty},\ \|\nn f\|_\infty\Big\},\\
 &\|\nn^2 P_{t,s} f\|_\infty\le c (s-t)^{-\ff 1 2} \|\nn f\|_\infty,\ \ 0\le t<t\le T.\end{align*}\end{cor}

 \beg{proof} By \eqref{BS} we have
 $$\|\nn P_{t,s}f\|_\infty\le c (t-s)^{-\ff 1 2} \|f\|_\infty$$ for some constant $c>0$. Next, by chain rule and \eqref{A*1},
 $$|\nn P_{t,s}f(x)|= \big|\E[\<\nn f(X_{t,s}^x), \nn X_{t,s}^x\>]\big|\le c \|\nn f \|_\infty,\ \ (t,s,x)\in D_T\times \R^d$$
holds  for some constant $c>0$. Moreover,
$$\nn P_{t,s}f(x)= \E[\<\nn f (X_{t,s}^x), \nn X_{t,s}^x\>]= \E[ g(X_{t,s}^x)],$$
where $g(X_{t,s}^x):= \big\<\nn f(X_{t,s}^x), \E(\nn X_{t,s}^x|X_{t,s}^x)\big\>.$ Combining this with \eqref{BS} and \eqref{A*1}, we find a constant $c>0$ such that
\beg{align*} &\|\nn^2 P_{t,s}f(x)\|\le \|\nn \E[g(X_{t,s}^x)]\|\\
&\le \ff 1 {s-t} \E\bigg[\big| g(X_{t,s}^x)\big|\cdot\bigg|\int_s^t \Big\<\big(\ss{2 a_{T-r}}\big)^{-1}(X_{t,r}^x) \nn_v X_{t,r}^x,\ \d W_r\Big\>\bigg|\bigg] \\
&\le \ff 1 {t-s} \big(\E|g(X_{t,s}^x)|^2\big)^{\ff 1 2} \bigg(\E\int_t^s  \|a^{-1}\|_\infty \|\nn X_{t,r}^x\|^2\d r\bigg)^{\ff 1 2} \le c \|\nn f\|_\infty.\end{align*}
Then the proof is finished.
 \end{proof}

\section{Proofs of Theorem \ref{T1}(2)-(3)}
We will need the following lemma implied by \cite[Theorem 2.1, Theorem 3.1, Lemma 3.3]{YZ}, see also \cite{XXZZ} and references within for the case $b^{(1)}=0$.

\beg{lem}\label{LYZ}  Assume $(A)(1),(A)(3)$ and   $\|b^{(0)}\|_{\tt L_{q_0}^{p_0}}<\infty$ for some
$(p_0,q_0)\in \scr K.$ Let $\si_t=\ss{2a_t}$. Then the following assertions hold.
\beg{enumerate} \item[$(1)$]
For any $p,q>1,$ $\ll\ge 0$, $0\le t_0<t_1\le T$ and $f\in \tt L_{q}^p(t_0,t_1)$,  the   PDE
\beq\label{1.2.6} (\pp_t +L_t)u_t=  \ll u_t+ f_t,\ \ t\in [t_0,t_1],  u_{t_1}=0,\end{equation}
  has a unique solution in $\tt H_q^{2,p}(t_0,t_1)$.
   If $(2p,2q)\in \scr K$, then
  there exist a constant  $  c>0$    such that for any $0\le t_0<t_1\le T$ and $f\in \tt L_{q}^p(t_0,t_1)$,
  the solution satisfies
\beg{align*}   \| u\|_{\infty}+ \|\nn u\|_\infty+ \|(\pp_t+\nn_{b^{(1)}}) u\|_{\tt L_q^p(t_0,t_1)}
 + \|\nn^2 u\|_{\tt L_q^p(t_0,t_1)} \le c \|f\|_{\tt L_q^p(t_0,t_1)}.\end{align*}
\item[$(2)$] Let   $(X_t)_{t\in [0,T]}$  be a continuous adapted process on $\R^d$ satisfying
\beq\label{1.2.1} X_t=X_0+\int_0^t b_s(X_s)\d s +\int_0^t \si_s  (X_s)\d W_s,\ \ t\in [0,T].\end{equation}
    For any $p,q>1$ with $(2p,2q)\in \scr K$,   there exists a constant $c>0$  such that for any $X_t$ satisfying $\eqref{1.2.1}$,
 \beg{align*} \E\bigg(\int_{t}^{s} |f_r(X_r)|\d r\bigg|\F_{t}\bigg)
 \le      c \|f\|_{\tt L_{q}^p(t,s)},\ \ (t,s)\in D_T, f\in \tt L_q^p(t,s).  \end{align*}
\item[$(3)$] Let $p,q> 1$ with $\ff d p+\ff 2 q<1$. For any $u\in\tt H_{q}^{2,p} $ with $\|(\pp_t +b^{(1)})u\|_{\tt L_q^p}<\infty$,
$\{u_t(X_t)\}_{t\in [0,T]}$ is a semimartingale satisfying
\beg{align*} \d u_t(X_t) =  L_t u_t(X_t)\d t + \big\<\nn u_t(X_t),\si_t(X_t)\d W_t\big\>,\ \ t\in [0,T].\end{align*}
\end{enumerate} \end{lem}

 In the following we    consider $E=\R^d$ and $\T^d$ respectively.

\subsection{$E=\R^d$}

\beg{proof}[Proof of Theorem \ref{T1}(2)]   Let $u\in \scr U(p_0,q_0)$ solve \eqref{NS}. Then
\beq\label{*D}  u\in \tt H_{q_0}^{2, p_0},\ \
  \|(\pp_t + b^{(1)}\cdot\nn) u\|_{\tt L_{q_0}^{p_0}}<\infty \end{equation}   as required by Lemma \ref{LYZ}(3).
It remains  to prove \eqref{SL}, which together with Corollary \ref{CN} implies \eqref{EST}.

Let
\beq\label{LS}\beg{split}  &  \scr L_t:= {\rm tr}\{a_{T-t}\nn^2\}+ \tt b_{t}\cdot\nn,\\
&\tt b_{t}(x):= b_{T-t}  (x)-\E u_0(X_{t,T}^x)-\E\int_{t}^T V_{T-s}(X_{t,s}^x)\d s,\ \ (t,x)\in [0,T]\times \R^d.\end{split}\end{equation}
Since $\|u_0\|_\infty+\int_0^T \|V_t\|_\infty\d t<\infty$,  $\|b^{(0)}\|_{\tt L_{q_0}^{p_0}}<\infty$ implies
   $\tt b_{t}(x):= b_{T-t}^{(1)}(x)+ \tt b_{t}^{(0)}  (x)$ with $\|\tt b^{(0)}\|_{\tt L_{q_0}^{p_0}}<\infty.$
 Then    $(A)$ holds for $\tt b$ replacing $b$, so that by \eqref{*D} and   Lemma \ref{LYZ}(3),   the following It\^o's formula holds for $X_{t,s}^x$ solving \eqref{S1}:
\beq\label{*1} \d u_{T-s}(X_{t,s}^x) = \big(\pp_s+\scr L_s\big) u_{T-s}(X_{t,s}^x) \d s +\big\{\nn u_{T-s}(X_{t,s}^x)\big\}^*\ss{2a_{T-s}(X_{t,s}^x)}\d W_s,\ \ s\in [t,T], \end{equation}
where $(\nn u)^*_{ij}:=(\pp_j u^i)_{1\le i,j\le d}.$
By \eqref{NS} and \eqref{LS},  we obtain
\beg{align*}&\big(\pp_s+\scr L_s\big) u_{T-s}(X_{t,s}^x)+V_{T-s}(X_{t,s}^x)\\
&=  \bigg\{\bigg[ u_{T-s}(y)-  \E u_0(X_{s,T}^y)- \E \int_{s}^TV_{T-r}(X_{s,r}^y)\d r\bigg]_{y=X_{t,s}^x}\cdot\nn\bigg\}
u_{T-s}(X_{t,s}^x).\end{align*}
Combining this with the follow property \eqref{FL} and
 \eqref{*1}, we derive
\beg{align*} &\E u_0(X_{t,T}^x) -  u_{T-t}(x) = \E\big[u_{T-T}(X_{t,T}^x)- u_{T-t}(X_{t,t}^x)\big]\\
&= \E\int_t^T \bigg\{\bigg( u_{T-s}(y)-  \E u_0(X_{s,T}^y) - \E\int_{s}^TV_{T-r}(X_{s,r}^y)\d r\bigg)_{y=X_{t,s}^x}\cdot\nn\bigg\}
u_{T-s}(X_{t,s}^x)\d s\\
&\quad - \E  \int_{t}^TV_{T-s}(X_{t,s}^x) \d s,\ \ \ (t,x)\in [0,T]\times \R^d.\end{align*}
Letting
$$h_t:= \sup_{x\in \R^d} \bigg|u_{T-t}(x)- \E u_0(X_{t,T}^x)-\E\int_{t}^T V_{T-s}(X_{t,s}^x)\d s\bigg|,\ \ t\in [0,T],$$
we arrive at
$$h_{t} \le \int_t^T h_{s} \|\nn u\|_\infty \d s,\ \ t\in [0,T].$$
By Grownwall's inequality we prove $h_t=0$ for $t\in [0,T],$ hence \eqref{SL} holds.\end{proof}

\beg{proof}[Proof of Theorem \ref{T1}(3)] (a) Let $P_{t,s}f= \E[f(X_{t,s}^x)]$ for $f\in \B_b(\R^d)$, where $X_{t,s}^x$ solves \eqref{S1}.
For $u$ given by \eqref{SL} we have
\beq\label{DC1} u_t= P_{T-t,T}u_0+\int_{T-t}^T P_{T-t,s} V_{T-s} \d s,\ \ t\in [0,T].\end{equation}
By $\|u_0\|_\infty+\int_0^T\|V_t\|_\infty\d t<\infty$ and \eqref{EST}, we find a constant $c>0$ such that
\beq\label{EST1}  \|u\|_\infty+\|\nn u\|_\infty\le c,\ \ \|\nn^2 u_t\|_\infty\le c t^{-\ff 1 2},\ \ t\in (0,T].\end{equation}
Moreover, the SDE \eqref{S1} becomes
\beq\label{S1'}  \beg{split}&\d X_{t,s}^x=\ss{2 a_{T-s}}(X_{t,s}^x)\d W_s
  +\big\{b_{T-s}-u_{T-s}\big\}(X_{t,s}^x) \d s,  \\
 &\qquad  \ t\in [0,T], s\in [t,T], X_{t,t}^x=x\in\R^d,\end{split}\end{equation}
and the generator in \eqref{LS} reduces to
$$\scr L_s:={\rm tr}\big\{a_{T-s}\nn^2\big\}+\big\{b_{T-s}-u_{T-s}\big\}\cdot\nn,\ \ s\in [0,T].$$

(b) We prove  the Kolmogorov backward  equation
\beq\label{DC2} \pp_t P_{t,s}f =-\scr L_t P_{t,s}f,\ \ f\in \C_b^2, t\in [0,s], s\in (0,T].\end{equation}
  For any $f\in \C_b^2$,
by It\^o's formula we have
\beq\label{KC} P_{t,s}f(x)= f(x)+\int_t^s P_{t,r}(\scr L_rf)(x)\d r,\ \ (t,s)\in D_T,\end{equation}
where $\int_t^s P_{t,r}(\scr L_r f)(x)\d r=\E \int_t^s \scr L_r f(X_{t,r}^x)\d r$ exists, since    Krylov's estimate in Lemma \ref{LYZ}(2) holds under $(A) $ and $\|u\|_\infty<\infty$.

By \eqref{KC}, we obtain  the Kolmogorov forward equation
\beq\label{KF} \pp_s P_{t,s}f= P_{t,s}(\scr L_s f),\ \  \ s\in [t,T].\end{equation}
On the other hand, $b^{(1)}=0$ and $(A)$ imply
\beq\label{*DD} \|\scr Lf\|_{\tt L_{q_0}^{p_0}}\le c_0\|f\|_{\C_b^2}\end{equation} for some constant $c_0>0$. By Lemma \ref{LYZ}(1),   for any $s\in (0,T]$, the PDE
\beq\label{PDE1} (\pp_t+\scr L_t)\tt u_t= - \scr L_t f,\ \  t\in [0,s], \tt u_s=0\end{equation}
has a unique solution $\tt u\in \scr U(p_0,q_0)$, such that
\beq\label{0B} \|\nn^2 \tt u\|_{\tt L_{q_0}^{p_0}(0,s)}\le c_1 \|\scr Lf\|_{\tt L_{q_0}^{p_0}(0,s)}
 \end{equation}
holds for some constant $c_1>0$  independent of $s$.
By  It\^o's formula in  Lemma \ref{LYZ}(3),
$$\d \tt u_t(X_{0,t}^x)= -\scr L_t f(X_{0,t}^x)+\big\<\nn f(X_{0,t}^x),\ss{2a_{T-t}}(X_{0,t}^x) \d W_t\big\>,\ \ t\in [0,s].$$
This and \eqref{KF} imply
\beg{align*} & 0 = \tt u_s(x)= \tt u_t(x) - \int_t^s (P_{t,r} \scr L_r f)(x)\d r \\
&=\tt u_t(x)- \int_t^s \ff{\d}{\d r} (P_{t,r}f)\d r= \tt u_t(x)- P_{t,s}f(x) + f(x),\ \ t\in [0,s].\end{align*}Thus,
\beq\label{1B} \tt u_t= P_{t,s} f-f,\ \ t\in [0,s].\end{equation}
Combining this with \eqref{PDE1} we derive \eqref{DC2}.

(c) By \eqref{EST1} and \eqref{DC2},   we see that $u$ solves \eqref{SL} with $u\in \scr U(p_0,q_0)$ provided
\beq\label{DV} \|\nn^2 u\|_{\tt L_{q_0}^{p_0}}<\infty.\end{equation}
By \eqref{*DD}, \eqref{0B} and \eqref{1B},
we find a constant $c_2>0$ such that
$$\sup_{t\in [0,s]}   \|  \nn^2 P_{\cdot,s} f \|_{\tt L_{q_0}^{p_0}(0,s)}
\le c_2  \|  f\|_{\C_b^2},\ \ s\in (0,T], f \in \C_b^2.$$
Combining this with \eqref{DC1}, $b^{(1)}=0$ and  $\|u_0\|_{\C_b^2}+\int_0^T\|V_t\|_{\C_b^2}\d t<\infty$,
we prove \eqref{DV}.

\end{proof}

\subsection{$E=\T^d$}

In this case, all functions on $E$ are extended to $\R^d$ as in \eqref{PR}, so that the proof for $E=\R^d$ works also for the present setting if we could verify the following periodic property for the solution of \eqref{S1}:
\beq\label{PP} X_{t,s}^{x+ k}= X_{t,s}^x+k,\ \ (t,s)\in D_T,\ x\in \R^d, \ k\in \Z^d.\end{equation}
Let $\tt X_{s,t}^x:= X_{t,s}^x+k.$ Since the coefficients of \eqref{S1} satisfies \eqref{PR}, $\tt X_{t,s}^x$
solves \eqref{S1} with $\tt X_{t,t}^x= x+k$. By the uniqueness of \eqref{S1} ensured by Theorem \ref{T1}(1), we derive \eqref{PP}.

\section{Application to \eqref{E1}}

 For any $n\in \mathbb N$, let $\C_b^n$ be the class of real functions $f$ on $E$ having derivatives up to order $n$ such that
$$\|f\|_{\C_b^n}:=\sum_{i=0}^n \|\nn^i f\|_\infty<\infty,$$
where $\nn^0f:=f$. Moreover, for $\aa\in (0,1)$, we denote $f\in \C_b^{n+\aa}$ if $f\in \C_b^n$ such that
$$\|f\|_{\C_b^{n+\aa}}:= \|f\|_{\C_b^n}+ \sup_{x\ne y} \ff{\|\nn^n f(x)-\nn^n f(y)\|}{|x-y|^\aa}<\infty.$$
 Consider  the following future distribution dependent SDE on $\R^d$:
\beq\label{SDE}  \d X_{t,s}^x= \bigg[\E\int_s^T\nn\wp_{T-r}(X_{s,r}^y)\d r -\E u_0(X_{s,T}^y)\bigg]_{y=X_{t,s}^x}\d s + \ss{2\kk}\d W_s,\ \ X_{t,t}^x=x, s\in [t,T].\end{equation}
See Definition \ref{D1} below for the definition of solution.  When $E=\T^d:=\R^d/\Z^d$, we extend $u_0$ and $\wp_t$ to $\R^d$ periodically, i.e. for a function $f$ on $\T^d$, it is extended to $\R^d$ as in \eqref{PR}.  
  With this extension, we also have  the SDE \eqref{SDE} for the case $E=\T^d$.

\beg{thm}\label{C2} If there exists $n\ge 2$ such that $u_0\in \C_b^{n}$ and $\wp_t\in \C_b^{n}$ for a.e. $t\in [0,T]$  with
$$\int_0^T \big(\| \nn \wp_t\|_\infty^2+\|\wp_t\|_{\C_b^n} \big)\d t <\infty.$$ Then $\eqref{SDE}$ is well-posed and  \eqref{E1} has a unique solution satisfying
\beq\label{EN} \sup_{t\in [0,T]} \|u_t\|_{\C_b^n}<\infty,\end{equation}
and the  solution is given by
\beq\label{*} u_t(x)=  \E u_0(X_{T-t,T}^x) -\E\int_{T-t}^T \nn \wp_{T-s}(X_{T-t,s}^x)\d s.\end{equation}

\end{thm}

 We only prove for $E=\R^d$ as the case for $E=\T^d$ follows by extending functions from $\T^d$  to $\R^d$ as in \eqref{PR}.

Let $I_d$ be the $d\times d$ identity matrix. By Theorem \ref{T1} with $b=0, a=\kk I_d$ and $V=-\nn\wp$, for any $(p_0,q_0)\in \scr K$, \eqref{E1} has a unique solution
in the class $\scr U(p_0,q_0)$, and by \eqref{*},
\beq\label{UUO}\beg{split}  u_t(x) &:= \E u_0(X_{T-t,T}^x) -\E\int_{T-t}^T \nn\wp_{T-s}(X_{T-t,s}^x)\d s\\
&= P_{T-t,T} u_0(x)-\int_{T-t}^T P_{T-t,s} \nn\wp_{T-s}(x)\d s, \ \ (t,x)\in [0,T]\times \R^d.\end{split}\end{equation} By \eqref{S1'} for the present $a$ and $b$, $X_{t,s}^x$ solves the SDE
\beq\label{SN} \d X_{t,s}^x= \ss{2\kk} \d W_s- u_{T-s}(X_{t,s}^x)\d s,\ \ X_{t,t}^x=x, t\in [0,T], s\in [t,T],\end{equation}
and the generator is
$$\scr L_s:=\kk \DD- u_{T-s}\cdot\nn,\ \ s\in [0,T].$$
It remains to prove \eqref{EN}.
To this end, we present the following lemma.

\beg{lem}\label{L} Let $P_{t,s}f:= \E [f(X_{t,s}^x)]$ for the SDE $\eqref{SN}$. Let
  $m\ge 1$ such that
\beq\label{LL)} \sup_{t\in [0,T]}\|u_t\|_{\C_b^m}+ \|f\|_{\C_b^{m+1}} <\infty,\end{equation}
  then $\sup_{(t,s)\in D_T} \|P_{t,s}f\|_{\C_b^{m+1}}<\infty.$ \end{lem}

  \beg{proof}   By \eqref{SN} and $ \sup_{t\in [0,T]}\|u_t\|_{\C_b^m}<\infty$, we have
  $$  \sup_{(t,s,x)\in D_T\times\R^d} \E\big[\|\nn^i X_{t,s}^x\|\big]<\infty,\ \ 1\le i\le m.$$
  By chain rule, this implies that for some constant $c_0>0$,
  \beq\label{JL} \sup_{(t,s)\in D_T} \|P_{t,s}g\|_{\C_b^m} \le c_0 \|g\|_{\C_b^m},\ \ g\in \C_b^m.\end{equation}
  Let $P_t^0=\e^{\kk\DD t}$. By $\pp_r P_{r-t}^0= P_{r-t}^0 \kk \DD$ and \eqref{DC2}, we have
  $$\pp_r P_{r-t}^0 P_{r,s}f= P_{r-t}^0\<\nn P_{r,s}f, u_{T-r}\>,\  \ r\in [t,s].$$ So,
  \beq\label{JL2} P_{t,s}f= P_{s-t}^0f- \int_t^s P_{r-t}^0\<\nn P_{r,s}f, u_{T-r}\>\d r.\end{equation}
  It is well known that for any $\aa,\bb\ge 0$ there exists a constant $c_{\aa,\bb}>0$ such that
  \beq\label{JL3} \|P_t^0 g\|_{\C_b^{\aa+\bb}}\le c_{\aa,\bb} t^{-\ff \aa 2} \|g\|_{\C_b^\bb},\ \ t>0, g\in \C_b^\bb.\end{equation}
  This together with \eqref{JL2} implies that for some constants $c_1,c_2>0$,
  $$\|P_{t,s}f\|_{\C_b^{m+\ff 1 2}}\le c_1 \|f\|_{\C_b^{m+\ff 1 2}} + c_1 \int_t^s (t+r-s)^{-\ff 3 4}
  \|\<\nn P_{r,s}f, u_{T-r}\>\|_{\C_b^{m-1}}\d r.$$
 Combining this  with \eqref{JL} and $\|f\|_{\C_b^m}+\sup_{t\in [0,T]}\|u_t\|_{\C_b^m}<\infty,$ we obtain
$$ \sup_{(t,s)\in D_T} \|P_{t,s}f\|_{\C_b^{m+\ff 1 2}}<\infty.$$ By this together with
   \eqref{JL2} and   \eqref{LL)}, we find a constant $c_2>0$ such that
\beg{align*}& \sup_{(t,s)\in D_T} \|P_{t,s}f\|_{\C_b^{m+1}}\le c_2 \|f\|_{\C_b^{m+1}} \\
&+ c_2 \sup_{(t,s)\in D_T} \int_t^s (t+r-s)^{-\ff 3 4}
  \|\<\nn P_{r,s}f, u_{T-r}\>\|_{\C_b^{m-\ff 12}}\d r<\infty.\end{align*}

  \end{proof}

  We now prove \eqref{EN} as follows. By $u\in \scr U(p_0,q_0),$ we have
  $$\|u\|_\infty+\|\nn u\|_\infty<\infty.$$
  Combining this with \eqref{UUO} and Lemma \ref{L}, we prove \eqref{EN} by inducing in $m$ up to $m=n$.

\end{document}